\documentclass[12pt, reqno]{amsart}
\usepackage{amsmath, amsthm, amscd, amsfonts, amssymb, graphicx, color}
\usepackage[bookmarksnumbered, colorlinks, plainpages]{hyperref}
\hypersetup{colorlinks=true,linkcolor=red, anchorcolor=green,
citecolor=cyan, urlcolor=red, filecolor=magenta, pdftoolbar=true}

\usepackage{amsmath}
\usepackage{amsfonts}
\usepackage{amssymb}
\usepackage{mathrsfs}
\usepackage{eufrak}

\newtheorem{theorem}{Theorem}[section]
\newtheorem{lemma}[theorem]{Lemma}
\newtheorem{proposition}[theorem]{Proposition}

\theoremstyle{definition}

\theoremstyle{remark}
\newtheorem{remark}[theorem]{Remark}
\numberwithin{equation}{section}

\theoremstyle{theorem}

\newtheorem{other}{Theorem}         

\newenvironment{pf}{\noindent{\emph{Proof.}}}{$\Box$}
\newenvironment{Pf}{\noindent{\emph{Proof of}}}{$\Box$}
\newcommand{\C}{{\mathbb C}}
\newcommand{\D}{{\mathbb D}}

\newcommand{\ig}{\stackrel{\text{def}}{=}}
\newcommand{\hol}{\mathcal Hol}
\newcommand{\Aut }{\text{Aut}}
\newcommand{\dbar}{\overline{\partial}}
 \DeclareMathOperator{\op}{o}

\begin{document}
\title[Multipliers and integration operators]{Multipliers and integration operators between conformally invariant
spaces}
\thanks{This research is supported in part by a grant from \lq\lq El Ministerio de
Ciencia, Innovaci\'{o}n y Universidades\rq\rq , Spain
(PGC2018-096166-B-I00) and by grants from la Junta de Andaluc\'{\i}a
(FQM-210 and UMA18-FEDERJA-002).}


\author{Daniel Girela}
\address{Departamento de An\'alisis Matem\'atico,
Estad\'{\i}stica e Investigaci\'{o}n Operativa, y Matem\'{a}tica
Aplicada, Universidad de M\'{a}laga, 29071 M\'{a}laga, Spain}
\email{girela@uma.es}
\author{Noel Merch\'{a}n}
\address{Departamento de Matem\'{a}tica Aplicada, Universidad de M\'{a}laga, 29071 M\'{a}laga, Spain}
\email{noel@uma.es}

\begin{abstract}
In this paper we are concerned with two classes of conformally
invariant spaces of analytic functions in the unit disc $\D$, the
Besov spaces $B^p$ $(1\le p<\infty )$ and the $Q_s$ spaces
$(0<s<\infty )$. Our main objective is to characterize for a given
pair $(X, Y)$ of spaces in these classes, the space of pointwise
multipliers $M(X, Y)$, as well as to study the related questions of
obtaining characterizations of those $g$ analytic in $\D $ such that
the Volterra operator $T_g$ or the companion operator $I_g$ with
symbol $g$ is a bounded operator from $X$ into $Y$. \end{abstract}


 \maketitle
\section{Introduction}\label{Intro}
\par Let $\D=\{z\in\C: |z|<1\}$ denote the open unit disc of the complex
plane $\C$ and let $\hol (\D)$ be the space of all analytic
functions in $\D$ endowed with the topology of uniform convergence
on compact subsets.
\par
If $\,0<r<1\,$ and $\,f\in \hol (\D)$, we set
\begin{align*}
M_p(r,f)\,=\,&\left(\frac{1}{2\pi }\int_0^{2\pi }
|f(re^{it})|^p\,dt\right)^{1/p}, \quad 0<p<\infty ,
\\
M_\infty(r,f)\,=\,&\sup_{\vert z\vert =r}|f(z)|.
\end{align*}
\par If $\,0<p\le \infty $\, the Hardy space $H^p$ consists of those
$f\in \hol(\mathbb D)$ such that \begin{equation*}\Vert f\Vert
_{H^p}\ig \sup_{0<r<1}M_p(r,f)<\infty .\end{equation*} We mention
\cite{D} for the theory of $H^p$-spaces.
\par If $0<p<\infty $ and $\alpha>-1$, the weighted
Bergman space $A^p_\alpha$ consists of those $f\in \hol(\mathbb D)$
such that
\[
\Vert f\Vert _{A^p_\alpha }\ig \left ((\alpha +1) \int_\mathbb
D(1-\vert z\vert )\sp\alpha \vert f(z)\vert ^p\, dA(z) \right
)^{1/p}<\infty .\] The unweighted Bergman space $A\sp p\sb 0 $ is
simply denoted by $A\sp p $. Here, $dA(z) =\frac{1}{\pi}dx\,dy $
denotes the normalized Lebesgue area measure in $\mathbb D$. We
refer to \cite{DS}, \cite{HKZ} and \cite{Zhu} for the theory of
these spaces.
\par\medskip
We let $\operatorname{Aut}(\D) $ denote the set of all disc
automorphisms, that is, of all one-to-one analytic maps $\varphi $
from $\D $ onto itself. It is well known that
$\operatorname{Aut}(\D) $ coincides with the set of all M\"obius
transformations from $\D $ onto itself:
$$ \operatorname{Aut}(\D )=\{ \lambda \varphi \sb a : \vert\lambda \vert =1, \,
a\in\D \, \} \,, $$ where $\varphi\sb a(z)=(a-z)/(1-\overline a z) $
($z\in \D$).
\par
A linear space $X $ of analytic functions in $\D $ is said to be
{\it conformally invariant\/} or {\it M\"obius invariant\/} if
whenever $f\in X $, then also $f\circ \varphi \in X $ for any
$\varphi\in\operatorname{Aut}(\D) $ and, moreover, $X $ is equipped
with a semi-norm $\rho $ for which there exists a positive constant
$C $ such that $$\rho (f\circ \varphi )\,\le \,C\rho (f ),
\quad\text{whenever $f\in X $ and $\varphi \in \Aut (\mathbb  D
)$}.$$ The articles \cite{AFP} and \cite{RT} are fundamental
references for the theory of M\"{o}bius invariant spaces which have
attracted much attention in recent years (see, e.\@\@g.,
\cite{AlSim,DGV1,DGV2,GM1,T,Zhu-MI,Zhu}).
\par
The {\it Bloch space\/} $\mathcal B $ consists of all analytic
functions $f $ in $\D $ such that $$\rho _{ _{\mathcal B}}(f)\ig
\sup\sb {z \in {\D }}\,(1-|z|\sp2)\,|f\sp\prime(z)|<\infty .$$ The
Schwarz-Pick lemma easily implies that $\rho _{ _{\mathcal B}}$ is a
conformally invariant seminorm, thus $\mathcal B$ is a conformally
invariant space. It is also a Banach space with the norm $\Vert
\cdot \Vert _{\mathcal B}$ defined by $\Vert f\Vert _{\mathcal
B}=\vert f(0)\vert \,+\,\rho _{ _{\mathcal B}}(f)$. The little Bloch
space $\mathcal B_0$ is the set of those $f\in \mathcal B$ such that
$\lim_{\vert z\vert \to 1}(1-\vert z\vert ^2)\vert f^\prime (z)\vert
=0$. Alternatively, $\mathcal B_0$ is the closure of the polynomials
in the Bloch norm. A classical reference for the theory of Bloch
functions is \cite{ACP}. Rubel and Timoney \cite{RT} proved that
$\mathcal B$ is the largest \lq\lq reasonable\rq\rq \, M\"{o}bius
invariant space. More precisely, they proved the following result.
\begin{other}\label{B-maximal} Let $X$ be a M\"{o}bius invariant
linear space of analytic functions in $\mathbb D $ and let $\rho $
be a M\"{o}bius invariant seminorm on $X $. If there exists a
non-zero decent linear functional $L $ on $X $ which is continuous
with respect to $\rho $, then $X\subset \mathcal B $ and there
exists a constant $A>0 $ such that $\rho _{ _{\mathcal
B}}(f)\,\le\,A\rho (f)$, for all $f\in X$.\end{other}
\par Here, a linear functional $L$ on $X$ is said to be decent if it
extends continuously to $\hol (\D )$.\par\medskip
  The space $BMOA $
consists of those functions $f $ in $H\sp1 $ whose boundary values
have bounded mean oscillation on the unit circle $\partial\D $ as
defined by F.~John and L.~Nirenberg. There are many
characterizations of $BMOA $ functions. Let us mention the
following:
\par\smallskip
{\it If  $f \in \hol (\D )$, then $f\in BMOA $ if and only if $
\Vert f\Vert \sb {BMOA}\ig\vert f(0)\vert \,+\,\rho_{*}(f)<\infty $,
where}
$$\rho_{*}(f)=\sup\sb
{a\in \D }\Vert f\circ \varphi \sb a-f(a)\Vert \sb {H\sp 2}.$$  It
is well known that $H\sp\infty \subset BMOA \subset\mathcal B $ and
that $BMOA$ equipped with the seminorm $\rho _*$ is a M\"{o}bius
invariant space. The space $VMOA$ consists of those $f\in BMOA$ such
that $ \lim _{\vert a\vert \to 1}\Vert f\circ \varphi \sb
a-f(a)\Vert \sb {H\sp 2}=0$, it is the closure of the polynomials in
the $BMOA$-norm. We mention \cite{G:BMOA} as a general reference for
the space $BMOA $.\par\medskip Other important M\"{o}bius invariant
spaces are the Besov spaces and the $Q_s$ spaces.
\par
For $1<p<\infty $, the {\it analytic Besov space\/} $B\sp p $ is
defined as the set of all functions $f $ analytic in $\mathbb D $
such that $f^\prime \in A^p_{p-2}$. All $B\sp p $ spaces
($1<p<\infty $) are conformally invariant with respect to the
semi-norm $\rho_ {_ {B^p}}$  defined by $$\rho_ {_ {B^p}}(f)\,\ig
\,\Vert f^\prime \Vert _{A^p_{p-2}}$$ (see \cite[p.\,\@112]{AFP} or
\cite[p.\,\@46]{DGV1}) and Banach spaces with the norm $\Vert \cdot
\Vert _{B^p}$ defined by  $\Vert f\Vert _{B^p}=\vert
f(0)\vert\,+\,\rho_ {_ {B^p}}(f)$. An important and well-studied
case is the classical {\it Dirichlet space\/} $B\sp 2 $ (often
denoted by $\mathcal D $) of analytic functions whose image has a
finite area, counting multiplicities.
\par
The space $B\sp1 $ requires a special definition: it is the space of
all analytic functions $f $ in $\D $ for which $f\sp{\prime \prime}
\in A\sp 1 $. Although the semi-norm $\rho$ defined by $\rho
(f)=\Vert f\sp{\prime \prime}\Vert _{A^1}$ is not conformally
invariant, the space itself is. An alternative definition of $B^1$
with a conformally invariant semi-norm is given in \cite{AFP}, where
it is also proved that $B^1$ is contained in any M\"{o}bius
invariant space. A lot of information on Besov spaces can be found
in \cite{AFP,DGV1,DGV2,HW1,Z,Zhu}. Let us recall that
\begin{align*}& VMOA \,\subsetneq \mathcal B_0,\quad BMOA\,\subsetneq
\mathcal B,\\ \quad B^1\,&\subsetneq B^p\,\subsetneq
\,B^q\,\subsetneq VMOA\,\subsetneq BMOA,\quad 1\,<\,p\,<\,q\,<\infty
.\end{align*}
\par
If $0\le s<\infty$, we say that $f\in Q_s$ if $f$ is analytic in
$\D$ and $$\sup\sb {a\in\D } \int\sb {\D }|f\sp\prime(z)|\sp2
g(z,a)\sp s\,dA(z)<\infty \,,$$ where $ g(z,a)=\log(|1-\overline{a}
z|/|a-z|) $ is the Green function of $\D $. These spaces were
introduced by Aulaskari and Lappan \cite{au-la94} while looking for
characterizations of Bloch functions (see \cite{X} for the case
$s=2$). For $s>1$, $Q_s$ is the Bloch space, $Q_1=BMOA$, and
\[\mathcal D\subsetneq Q_{s_1}\subsetneq Q_{s_2}\subsetneq BMOA,\qquad0<s_1<s_2<1.\]
It is well known \cite{AXZ,Str} that for every $s$ with $0\le
s<\infty $, a function $f\in \hol (\D )$ belongs to $Q_s$ if and
only if
\begin{equation*}\rho_{ _{Q_s}}(f)\,\ig \,\left (
\sup_{a\in\D}\int_\D|f'(z)|^2(1-\vert \varphi _a(z)\vert
^2)^s\,dA(z)\right )^{1/2}\,<\,\infty.\end{equation*} All $Q_s$
spaces ($0\le s<\infty $) are conformally invariant with respect to
the semi-norm $\rho_{ _{Q_s}}$. They are also Banach spaces with the
norm $\Vert \cdot \Vert _{Q_s}$ defined by $\Vert f\Vert
_{Q_s}=\vert f(0)\vert \,+\,\rho_{ _{Q_s}}(f)$. We mention
\cite{X2,X3} as excellent references for the theory of $Q_s$-spaces.
\par Let us recall the following two facts which were first observed
in \cite{AC}.
\begin{equation}\label{BpQsple2} \text{If $0<p\le 2$,\, then\,
$B^p\subset Q_s$\, for all $s>0$.}\end{equation}
\begin{equation}\label{BpQsp>2} \text{If $2<p<\infty $,\, then\,
$B^p\subset Q_s$\, if and only if \,
$1-\frac{2}{p}<s$.}\end{equation}

\par\medskip For $g$ analytic in $\D$, the Volterra operator
$T_g$ is defined as follows:
\[
T_g(f)(z)\ig\int_0^zg'(\xi)f(\xi)\,d\xi, \,\,\, f\in \hol (\D
),\,\,\, z\in \D.\] We define also the companion operator $I_g$ by
\[
I_g(f)(z)\ig\int_0^zg(\xi)f'(\xi)d\xi,\,\,\, f\in \hol (\D ),\,\,\,
z\in \D.\] The integration operators $T_g$ and $I_g$ have been
studied in a good number of papers. Let us just mention here that
Pommerenke \cite{Pom} proved that $T_g$ is bounded on $H^2$ if and
only if $g\in BMOA$ and that Aleman and Siskakis \cite{AlSis-CVTA}
characterized those $g\in \hol (\D )$ for which $T_g$ is bounded on
$H^p$ ($p\ge 1$), while
 Aleman and Cima characterized in \cite{AlCi} those
$g\in \hol (\D )$ for which $T_g$ maps $H^p$ into $H^q$. Aleman and
Siskakis \cite{AlSis-Indiana} studied the operators $I_g$ and $T_g$
acting on Bergman spaces.
\par For $g\in\hol (\D )$, the
multiplication operator $M_g$ is defined by \[ M_g(f)(z)\ig
g(z)f(z),\quad f\in \hol (\D),\,\, z\in \D.
\]
If $X$ and $Y$ are two Banach spaces of analytic function in $\D $
continuously embedded in $\hol (\D )$ and $g\in\hol (\D )$ then $g$
is said to be a multiplier from $X$ to $Y$ if $M_g(X)\subset Y$. The
space of all multipliers from $X$ to $Y$ will be denoted by $M(X,Y)$
and $M(X)$ will stand for $M(X,X)$. Using the closed graph theorem
we see that for the three operators $T_g$, $I_g$, $M_g$, we have
that if one of them maps $X$ into $Y$ then it is continuous from $X$
to $Y$. We remark also that
\begin{equation}\label{IJM}
T_g(f)+I_g(f)=M_g(f)-f(0)g(0).\end{equation} Thus if two of the
operators $T_g, I_g, M_g$ are bounded from $X$ to $Y$ so is the
third one.
\par\medskip
It is well known that if $X$ is nontrivial then $M(X)\subset
H^\infty $ (see, e.\,\@g., \cite[Lemma\,\@1.\,\@1]{ADMV} or
\cite[Lemma\,\@1.\,\@10]{Vi}), but $M(X, Y)$ need not be included in
$H^\infty $ if $Y\not\subset X$. However, when dealing with
M\"{o}bius invariant spaces we have the following result.
\begin{proposition}\label{MXY-Conf} Let $X$ and $Y$ be two
M\"{o}bius invariant spaces of analytic functions in $\D $ equipped
with the seminorms $\rho _{ _X}$ and $\rho _{ _Y}$, respectively.
Suppose that there exists a non-trivial decent linear functional $L$
on $Y$ which is continuous with respect to $\rho _{ _Y}$. Let $g\in
\hol (\D )$. Then the following statements hold.
\begin{itemize} \item [(i)] If\, $M_g$ is continuous from $(X, \rho _{ _X})$ into $(Y, \rho
_{ _Y})$, then $g\in H^\infty $.
\item [(ii)] If\, $I_g$ is continuous from $(X, \rho _{ _X})$ into $(Y, \rho
_{ _Y})$, then $g\in H^\infty $.
\end{itemize}
\end{proposition}
\par\medskip Before embarking into the proof of
Proposition\,\@\ref{MXY-Conf}, let us mention that, as usual,
throughout the paper we shall be using the convention that $C=C(p,
\alpha ,q,\beta , \dots )$ will denote a positive constant which
depends only upon the displayed parameters $p, \alpha , q, \beta
\dots $ (which sometimes will be omitted) but not  necessarily the
same at different occurrences. Moreover, for two real-valued
functions $E_1, E_2$ we write $E_1\lesssim E_2$, or $E_1\gtrsim
E_2$, if there exists a positive constant $C$ independent of the
arguments such that $E_1\leq C E_2$, respectively $E_1\ge C E_2$. If
we have $E_1\lesssim E_2$ and  $E_1\gtrsim E_2$ simultaneously then
we say that $E_1$ and $E_2$ are equivalent and we write $E_1\asymp
E_2$. Also, if \,$1<p<\infty $, $p^\prime $ will stand for its
conjugate exponent, that is, $\frac{1}{p}+\frac{1}{p^\prime }=1$.
\par\medskip
\begin{Pf}{\it Proposition\,\@\ref{MXY-Conf}.} Since $X$ is conformally invariant, $\Aut (\mathbb D)\subset X $ \cite[p.\,\@114]{AFP} and
\begin{equation}\label{rhophia}\rho _{ _X}(\varphi_a)\asymp 1,\quad
a\in \mathbb D.\end{equation}
\par Suppose that $M_g$ is continuous from $(X, \rho _{ _X})$ into $(Y, \rho
_{ _Y})$. Using this, Theorem\,\@\ref{B-maximal}, and
(\ref{rhophia}) we obtain
\begin{equation*}\rho_{ _\mathcal B}(g\,\varphi _a)\,\lesssim \,\rho_{ _Y}
(g\,\varphi _a)\,\lesssim\,\rho _{ _X}(\varphi_a)\,\lesssim 1,\quad
a\in \mathbb D.
\end{equation*}
This implies that
$$(1-\vert a\vert ^2)\left \vert g^\prime (a)\varphi
_a(a)\,+\,g(a)\varphi_a^\prime (a)\right \vert \,\lesssim 1,\quad
a\in \mathbb D.$$ Since $\varphi (a)=0$ and $\varphi _a^\prime
(a)=-(1-\vert a\vert ^2)^{-1}$, it follows that $$\vert g(a)\vert
\,\lesssim 1,\quad a\in \mathbb D,$$ that is, $g\in H^\infty $.
\par Similarly, if we assume that $I_g$ is continuous from $(X, \rho _{ _X})$ into $(Y, \rho
_{ _Y})$, we obtain
$$\rho_{ _\mathcal B}\left (I_g(\varphi _a)\right )\,\lesssim
1,\quad a\in \mathbb D.$$ This implies that $$(1-\vert a\vert
^2)\left \vert \left (I_g(\varphi _a)\right )^\prime (a)\right \vert
\,=\,(1-\vert a\vert ^2)\vert \varphi_a^\prime (a)\vert \vert
g(a)\vert \,=\,\vert g(a)\vert \,\lesssim 1,\quad a\in \mathbb  D.$$
\end{Pf}
\par\medskip For notational convenience, set
$$\mathcal B\mathcal Q\,=\,\{ Q_s : 0\,\le \,s\,<\,\infty \}\,\cup \{ B^p :
1\,\le\,p\,<\infty \} .$$ The main purpose of this paper is
characterizing, for a given pair of spaces $X, Y\in \mathcal
B\mathcal Q$, the functions $g\in \hol (\D )$ such that the
operators $M_g$, $T_g$ and/or $I_g$ map $X$ into $Y$. When $X$ and
$Y$ are Besov spaces this question has been extensively studied
(see, e.\,\@g. \cite{ARS,GaGiPe-TAMS-2011,GP-JFA-06,Steg,Wu,Zo}).
Thus we shall concentrate ourselves to study these operators when
acting between a certain Besov space $B^p$ and a certain $Q_s$ space
and when acting between
 $Q_{s_1}$ and
$Q_{s_2}$ for a certain pair of positive numbers $s_1, s_2$.

\par\medskip
\section{Multipliers and integration operators from Besov spaces
into $Q_s$-spaces}\label{MIOBpQs}
\par For $\alpha >0$, {\it the $\alpha $-logarithmic Bloch space}
 $\mathcal
B_{\log ,\alpha }$ is the Banach space of those functions $f\in\hol
(\D )$ which satisfy
\begin{equation}\label{logalpha}\Vert f\Vert _{\log ,\alpha }\ig
\vert f(0)\vert +\sup_{z\in \mathbb D}(1-\vert z\vert ^2)\left (\log
\frac{2}{1-\vert z\vert ^2}\right )^\alpha \vert f^\prime (z)\vert
<\infty .\end{equation} For simplicity, the space $\mathcal B_{\log
, 1}$ will be denoted by $\mathcal B_{\log }$.\par
 It is clear that $B_{\log ,\alpha }\subset
\mathcal B_0$, for all  $\alpha >0$. Using the characterization of
$VMOA$ in terms of Carleson measures \cite[p.\,\@102]{G:BMOA}, it
follows easily that $$B_{\log ,\alpha }\subset VMOA, \quad\text{for
all $\alpha >1/2$}.$$ In particular, $\mathcal B_{\log }\subset
VMOA$.
\par  Brown and Shields
\cite{BS} showed that $M(\mathcal B)=\mathcal B_{\log }\cap H^\infty
$. The spaces $M(B^p, \mathcal B)$ ($1\le p<\infty $) were
characterized in \cite{GaGiMa-CAOT}. Namely, Theorem\,\@\,\@1 of
\cite{GaGiMa-CAOT} asserts that $M(B^1, \mathcal B)=H^\infty $ and
\begin{equation}\label{BpB}M(B^p, \mathcal B)=H^\infty \cap \mathcal
B_{\log ,1/p^\prime},\quad 1\,<\,p\,<\,\infty ,\end{equation} where
$p^\prime $ is the exponent conjugate to $p$, that is,
$\frac{1}{p}+\frac{1}{p^\prime }=1$.
\par In this section we extend
these results. In particular, we shall obtain for any pair $(p, s)$
with $2\,<\,p\,<\,\infty $ and $0\,<\,s\,<\,\infty $ a complete
characterization of the space of multipliers $M(B^p, Q_s)$.

\par\medskip Let us start with the case $s\ge 1$ which is the simplest one.
\begin{theorem}\label{BpBloch}
Let $g\in \hol (\D )$. Then:
\begin{itemize}
\item[(i)] $I_g$ maps $B^1$ into $\mathcal B$ if and only if $g\in
H^\infty $.
\item[(ii)] $M_g$ maps $B^1$ into $\mathcal B$ if and only if $g\in
H^\infty $.
\item[(iii)] $T_g$ maps $B^1$ into $\mathcal B$ if and only if $g\in
\mathcal B$.
\end{itemize}
\end{theorem}
\begin{pf} If $I_g(B^1)\subset \mathcal B$ then, using
Proposition\,\@\ref{MXY-Conf}, it follows that $g\in H^\infty $.
\par To prove the converse it suffices to recall that $B^1\subset \mathcal B$. Indeed, suppose that $g\in H^\infty $ and take
$f\in B^1$. Then
$$(1-\vert z\vert^2)\left\vert \left (I_g(f)\right )^\prime
(z)\right \vert \,=(1-\vert z\vert^2)\vert f^\prime (z)\vert \vert
g(z)\vert \,\le \,\Vert f\Vert _{\mathcal B}\Vert g\Vert _{H^\infty
}.$$ Thus $I_g(f)\in \mathcal B$.
\par Hence (i) is proved. Now, (ii) is contained in
\cite[Theorem\,\@\,\@1]{GaGiMa-CAOT}.
\par It remains to prove (iii). If $T_g(B^1)\subset \mathcal B$ then
$T_g(1)=g\,-\,g(0)\,\in \,\mathcal B$ and, hence $g\in \mathcal B$.
Conversely, if $g\in \mathcal B$ and $f\in \mathcal B^1$ then, using
the fact that $B^1\subset H^\infty $, we obtain
$$(1-\vert z\vert ^2)\left\vert \left (T_g(f)\right )^\prime
(z)\right \vert \,=(1-\vert z\vert^2)\vert g^\prime (z)\vert \vert
f(z)\vert \,\le \,\Vert g\Vert _{\mathcal B}\Vert f\Vert_{H^\infty
}.$$ Thus $T_g(f)\in \mathcal B$. Hence (iii) is also proved.
\end{pf}
\par\medskip
\begin{theorem}\label{Bpp>1Bloch}
Suppose that $1<p<\infty $, $\frac{1}{p}+\frac{1}{p^\prime }=1$, and
let $g\in \hol (\D )$. Then:
\begin{itemize}
\item[(i)] $I_g$ maps $B^p$ into $\mathcal B$ if and only if $g\in
H^\infty $.
\item[(ii)] $M_g$ maps $B^p$ into $\mathcal B$ if and only if $g\in
H^\infty\cap \mathcal B_{\log ,1/p^\prime }$.
\item[(iii)] $T_g$ maps $B^p$ into $\mathcal B$ if and only if $g\in
\mathcal B_{\log , 1/p^\prime }$.
\end{itemize}
\end{theorem}
\begin{pf}
If $I_g$ maps $B^p$ into $\mathcal B$ then
Proposition\,\@\ref{MXY-Conf} implies that $g\in H^\infty $.
Conversely, using that $B^p\subset \mathcal B$, we see that if $g\in
H^\infty $ and $f\in B^p$ then
$$(1-\vert z\vert ^2)\left \vert \left (I_g(f)\right )^\prime
(z)\right \vert \,=\,(1-\vert z\vert ^2)\vert f^\prime (z)\vert
\vert g(z)\vert \,\le \,\Vert f\Vert _{\mathcal B}\Vert g\Vert
_{H^\infty }.$$ Hence, $I_g(f)\in \mathcal B$. Thus (i) is proved
and (ii) reduces to (\ref{BpB}).
\par Finally, (iii) follows from the following more precise result.
\begin{theorem}\label{TgBpBloch} Suppose that $1<p<\infty $,
$\frac{1}{p}+\frac{1}{p^\prime }=1$, and let $g\in \hol (\D )$. Then
the following conditions are equivalent.
\begin{itemize}
\item[(a)] $T_g$ maps $B^p$ into $\mathcal B$.
\item[(b)] $g\in \mathcal B_{\log , 1/p^\prime }$.
\item[(c)] $T_g$ maps $B^p$ into $\mathcal B_0$.
\end{itemize}
\end{theorem}
{\it Proof of Theorem\,\@\ref{TgBpBloch}.} (a)\,$\Rightarrow $\,(b)
Suppose (a). By the closed graph theorem $T_g$ is a bounded operator
from $B^p$ into $\mathcal B$, hence
\begin{equation}\label{fff}(1-\vert z\vert ^2)\vert g^\prime
(z)f(z)\vert \lesssim \Vert f\Vert _{B^p}, \quad z\in \D,\,\,f\in
B^p.\end{equation} For $a\in \mathbb D$ with $a\neq 0$, set
\begin{equation}\label{faBp}f_a(z)\,=\,\left (\log\frac{1}{1-\vert a\vert ^2}\right
)^{-1/p}\log\frac{1}{1-\overline a\,z},\quad z\in \mathbb
D.\end{equation} It is readily seen that $f_a\in B^p$ for all $a$
and that $\Vert f_a\Vert _{B^p}\asymp 1$. Using this and taking
$f=f_a$ and $z=a$ in (\ref{fff}), we obtain
$$(1-\vert a\vert ^2)\vert g^\prime (a)\vert \left (\frac{1}{1-\vert
a\vert ^2)}\right )^{1/p^\prime }\,\lesssim 1,$$ that is $g\in
\mathcal B_{\log , 1/p^\prime }$.
\par (b)\,$\Rightarrow $\,(c) Suppose (b) and take $f\in B^p$. It
is well known that  $$\vert f(z)\vert \,=\,\op \left (\left (\log
\frac{1}{1-\vert z\vert ^2}\right )^{1/p^\prime }\right
),\quad\text{as $\vert z\vert \to 1$},$$ (see, e.\,\@g.,
\cite{HW1,Z}). This and (b) immediately yield that $T_g(f)\in
\mathcal B_0$.
\par The implication (c)\,$\Rightarrow $\,(a) is trivial. Hence the
proof of Theorem\,\@\ref{TgBpBloch} is finished and, consequently,
Theorem\,\@\ref{Bpp>1Bloch} is also proved.
\end{pf}
\par\medskip Let us turn now to the case $0<s\le 1$.
We shall consider first the Volterra operators $T_g$. For
$0<s<\infty $ and $\alpha
>0$ we set $$Q_{s, \log , \alpha }\,=\,\left \{ f\in \hol (\D ) :
\sup _{a\in \D }\left (\log\frac{2}{1-\vert a\vert }\right
)^{2\alpha }\int_{\D}\vert f^\prime (z)\vert ^2(1-\vert\varphi
_a(z)\vert ^2)^s\,dA(z)<\infty \right \}.$$ We have the following
results.
\begin{theorem}\label{TgBpQs} Suppose that $0<s\le 1$ and let $g\in
\hol (\D )$. Then:
\begin{itemize}
\item[(i)] $T_g$ maps $B^1$ into $Q_s$ if and only if $g\in Q_s$.
\item[(ii)] If $1<p<\infty $, $0<s\le 1$, and $T_g$ maps $B^p$ into
$Q_s$, then $g\in Q_{s, \log , 1/p^\prime }$.
\item[(iii)] If $1<p<\infty $, then $T_g$ maps $B^p$ into $Q_1=BMOA$ if and only if
$g\in Q_{1, \log , 1/p^\prime }$.
\item[(iv)] If $2<p<\infty $, $0<s<1$,  and $1-\frac{2}{p}<s$ then
$T_g$ maps $B^p$ into $Q_s$ if and only if $g\in Q_{s, \log ,
1/p^\prime }$.
\end{itemize}
\end{theorem}
\par\medskip
Before we get into the proofs of these results we shall introduce
some notation and recall some results which will be needed in our
work.
\par
\par If  $I\subset \partial\D$ is an
interval, $\vert I\vert $ will denote the length of $I$. The
\emph{Carleson square} $S(I)$ is defined as
$S(I)=\{re^{it}:\,e^{it}\in I,\quad 1-\frac{|I|}{2\pi }\le r <1\}$.
Also, for $a\in \D$, the Carleson box $S(a)$ is defined by
\begin{displaymath}
S(a)=\Big \{ z\in \D : 1-|z|\leq 1-|a|,\, \Big |\frac{\arg
(a\bar{z})}{2\pi}\Big |\leq \frac{1-|a|}{2} \Big \}.
\end{displaymath}
\par If $\, s>0$ and $\mu$ is a positive Borel  measure on  $\D$,
we shall say that $\mu $
 is an $s$-Carleson measure
  if
there exists a positive constant $C$ such that
\[
\mu\left(S(I)\right )\le C{|I|^s}, \quad\hbox{for any interval
$I\subset\partial\D $},
\]
or, equivalently,
 if there exists $C>0$ such that
\[
\mu\left(S(a)\right )\le C{(1-\vert a\vert )^s}, \quad\hbox{for all
$a\in\D$}.
\]
A $1$-Carleson measure will be simply called a Carleson measure.
\par\medskip These concepts were generalized in \cite{Zhao} as
follows: If $\mu$ is a positive Borel measure in $\D$, $0\le \alpha
<\infty $, and $0<s<\infty $, we say that $\mu$ is an
 $\alpha$-logarithmic $s$-Carleson measure
  if there exists a positive
 constant $C$ such that
 \[\frac{
\mu\left(S(I)\right )\left(\log \frac{2\pi }{\vert I\vert }\right
)^\alpha }{|I|^s}\le C, \quad\hbox{for any interval
$I\subset\partial\D $}
\] or, equivalently,
 if
$$
\sup_{a\in\D}\frac{\mu\left(S(a)\right)\left (\log
\frac{2}{1-|a|^2}\right )^{\alpha}}{(1-|a|^2)^s}<\infty .$$ Carleson
measures and logarithmic Carleson measures are known to play a basic
role in the study of the boundedness of a great number of operators
between analytic function spaces. In particular we recall the
Carleson embedding theorem for Hardy spaces which asserts that if
$0<p<\infty $ and $\mu$ is a positive Borel measure on $\D$ then
$\mu $ is a Carleson measure if and only if the Hardy space $H^p$ is
continuously embedded in $L^p(d\mu )$ (see \cite[Chapter\,\@9]{D}).
\par
In the next theorem we collect a number of known results which will
be needed in our work.
\begin{other}\label{Collect-Results}\par
\begin{itemize}
\item[(i)] If $0<s\le 1$ and $f\in \hol (\D )$, then $f\in Q_s$ if
and only if the Borel measure $\mu $ on $D$ defined by
$$d\mu (z)\,=\, (1-\vert z\vert ^2)^s\vert f^\prime (z)\vert
^2\,dA(z)$$ is an $s$-Carleson measure.
\item[(ii)] If $0\le\alpha <\infty $, $0<s<\infty $, and $\mu $ is a
positive Borel measure on $\D $ then $\mu $ is an $\alpha
$-logarithmic $s$-Carleson measure if and only if
$$\sup _{a\in \D }\left (\log \frac{2}{1-\vert a\vert ^2}\right
)^{\alpha }\int_{\D }\left (\frac{1-\vert a\vert ^2}{\vert
1-\overline a\,z\vert ^2}\right )^s\,d\mu (z)\,<\,\infty .$$
\item[(iii)] If $1<p\le 2$ then $B^p\subset Q_s$ for all $s>0$.
\item[(iv)] If $2<p<\infty $ and $1-\frac{2}{p}<s$, then
$B^p\,\subset \,Q_s$.
\item[(v)]
For $s>-1$, we let $\mathcal D_s$ be the space of those functions
$f\in \hol (\D )$ for which
$$\Vert f\Vert _{\mathcal D_s}\,\ig \,\vert f(0)\vert \,+\,\left
(\int _{\D }(1-\vert z\vert ^2)^s\vert f^\prime (z)\vert
^2\,dA(z)\right )^{1/2}\,<\,\infty .$$ \par Suppose that $0<s<1$ and
$\alpha
>1$, and let $\mu $ be a positive Borel measure on $\D $. If $\mu $
is an $\alpha $-logarithmic $s$ Carleson measure, then $\mu $ is a
Carleson measures for $\mathcal D_s$, that is, $\mathcal D_s$ is
continuously embedded in $L^2(d\mu )$.
\end{itemize}
\end{other}
\par Let us mention that (i) is due to Aulaskari, Stegenga and Xiao
\cite{ASX}, (ii) is due to Zhao \cite{Zhao}, (iii) and (iv) were
proved by Aulaskari and Csordas in \cite{AC}, and (v) is due to Pau
and Pel\'{a}ez \cite[Lemma\,\@1]{Pau-Pe-MZ-2009}.

Using Theorem\,\@\ref{Collect-Results}\,\@(ii) and the fact that
$$1-\vert \varphi _(z)\vert ^2\,=\,\frac{(1-\vert a\vert ^2)(1-\vert z\vert
^2)}{\vert 1-\overline a\,z\vert ^2},$$ we see that for a function
$f\in \hol (\D )$ we have  that $f\in Q_{s, \log , \alpha }$ if and
only if the measure $\mu $ defined by $d\mu (z)=(1-\vert z\vert
^2)^s\vert f^\prime (z)\vert ^2\,dA(z)$ is a $2\alpha $-logarithmic
$s$-Carleson measure. \par\medskip

\begin{Pf}{\it Theorem\,\@\ref{TgBpQs}\,\@(i).} Suppose that $T_g$ maps $B^1$ into $Q_s$. Since the
constant functions belong to $B^1$, we have that
$T_g(1)\,=\,g\,-\,g(0)\,\in \,Q_s$ and, hence, $g\in Q_s$. \par To
prove the converse, suppose that $g\in Q_s$. Then the measure $\mu $
defined by $$d\mu (z)\,=\,(1-\vert z\vert ^2)^s\vert g^\prime
(z)\vert ^2\,dA(z)$$ is an $s$-Carleson measure. Take now $f\in
B^1$, then $f\in H^\infty $ and, hence,
$$(1-\vert z\vert^2)^s\left \vert \left (T_g(f)\right )^\prime
(z)\right \vert ^2\,=\,(1-\vert z\vert^2)^s\vert g^\prime (z)\vert
^2\vert f(z)\vert ^2\,\le\, \Vert f\Vert_{H^\infty }^2(1-\vert
z\vert^2)^s\vert g^\prime (z)\vert ^2.$$ Since $\mu $ is an
$s$-Carleson measure, it follows readily that the measure $\nu $
given by $d\nu (z)\,=\, (1-\vert z\vert^2)^s\left \vert \left
(T_g(f)\right )^\prime (z)\right \vert ^2\,dA(z)$ is also an
$s$-Carleson measure and, hence, $T_g(f)\in Q_s$.
\end{Pf}
\par\medskip
\begin{Pf}{\it Theorem\,\@\ref{TgBpQs}\,\@(ii).}
\par
Suppose that $0<s\le 1$, $1<p<\infty $, and that $T_g$ maps $B^p$
into $Q_s$. For $a\in \D \setminus \{ 0\} $, set $$f_a(z)\,=\,\left
(\log \frac{1}{1-\vert a\vert ^2}\right
)^{-1/p}\log\frac{1}{1-\overline a\,z},\quad z\in \D,$$ as in
(\ref{faBp}). We have that $\Vert f_a\Vert _{B^p}\asymp 1$ and it is
also clear that
$$\vert f_a(z)\vert \,\asymp \,\left (\log \frac{1}{1-\vert a\vert
^2}\right )^{1/p^\prime },\quad z\in S(a).$$ Using these facts, we
obtain
\begin{align*}&\frac{\left (\log
\frac{1}{1-\vert a\vert ^2}\right )^{2/p^\prime }}{(1-\vert
a\vert^2)^s}\int_{S(a)}(1-\vert z\vert ^2)^s\vert g^\prime (z)\vert
^2\,dA(z)\\ \asymp &\,\frac{1}{(1-\vert
a\vert^2)^s}\int_{S(a)}(1-\vert z\vert ^2)^s\vert g^\prime
(z)f_a(z)\vert ^2\,dA(z)\\ = &\, \frac{1}{(1-\vert
a\vert^2)^s}\int_{S(a)}(1-\vert z\vert ^2)^s\vert \left
(T_g(f_a)\right )^\prime (z)\vert ^2\,dA(z).
\end{align*}
The fact that $T_g$ is a bounded operator from $B^p$ into $Q_s$,
implies that the measures $(1-\vert z\vert ^2)^s\vert\left
(T_g(f_a)\right )^\prime (z)\vert ^2\,dA(z)$ are $s$-Carleson
measures with constants controlled by $\Vert T_g\Vert ^2$. Then it
follows that the measure $(1-\vert z\vert ^2)^s\vert g^\prime
(z)\vert ^2\,dA(z)$ is a $2/p^\prime $-logarithmic $s$-Carleson
measure and, hence, $g\in Q_{s, \log , 1/p^\prime }.$
\end{Pf}
\par\medskip
\begin{Pf}{\it Theorem\,\@\ref{TgBpQs}\,\@(iii) and (iv).}
 In view of (ii) we only have to prove that if $g\in Q_{s, \log , 1/p^\prime
 }$ then $T_g$ maps $B^p$ into $Q_s$. \par Hence, take
$g\in Q_{s, \log , 1/p^\prime }$  and set $$K(g)= \sup _{a\in \D
}\left (\log\frac{2}{1-\vert a\vert }\right )^{2/p^\prime
}\int_{\D}\vert g^\prime (z)\vert ^2(1-\vert\varphi _a(z)\vert
^2)^s\,dA(z),$$ and  take $f\in B^p$. Set $F=T_g(f)$, we have to
prove that $F\in Q_s$ or, equivalently, that the measure $\mu_{
_{F}}$ defined by
$$d\mu_{ _{F}}(z)\,=\,(1-|z|^2)^s\vert F^\prime(z)|^2\,dA(z)$$ is an $s$-Carleson measure.
Let $a\in \D$. Using the well known fact that
$$1-\vert a\vert ^2\,\asymp \,\vert
1-\overline a\,z\vert ,\quad z\in S(a),$$ we obtain

\begin{align}\label{muFsCar}
&\frac{1}{(1-\vert a\vert ^2)^s}\int_{S(a)}
|F^\prime(z)|^2(1-|z|^2)^s\,dA(z)\,\asymp \,\int_{S(a)}
|F^\prime(z)|^2\frac{(1-|z|^2)^s(1-\vert a\vert^2)^s}{\vert
1-\overline a\,z\vert ^{2s}}\,dA(z)\nonumber \\ & = \,\,
\int_{S(a)} |f(z)|^2 |g^\prime(z)|^2 (1-\vert \varphi _a(z)\vert^2)^s\,dA(z) \nonumber\\
& \le \,\, 2\,\int_{\D}|f(a)|^2 |g^\prime(z)|^2(1-\vert \varphi
_a(z)\vert^2)^s\,dA(z) \nonumber\\ & \qquad +\,\, 2\,\int_{\D
}|f(z)-f(a)|^2 |g^\prime(z)|^2 (1-\vert \varphi
_a(z)\vert^2)^s\,dA(z) \nonumber \\ & = \,\,  2T_1(a)\,+\,2T_2(a).
\end{align}
Using the fact that \begin{equation}\label{grBp}\vert f(a)-f(0)\vert
\lesssim \Vert f\Vert _{B^p}\left (\log \frac{2}{1-\vert a\vert
^2}\right )^{1/p^\prime },\end{equation} we obtain
\begin{align}\label{T1BMOA}T_1(a)
 \lesssim & \,\, \Vert f\Vert
_{B^p}^2\left (\log \frac{2}{1-\vert a\vert^2}\right
 )^{2/p^\prime }\int_{\D }\vert g^\prime (z)\vert ^2(1-\vert
\varphi _a(z)\vert ^2)^s\,dA(z)  \lesssim \,\,  K(g)\Vert f\Vert
_{B^p}^2.
\end{align}
\par
To estimate $T_2(a)$ we shall treat separately the cases $s=1$ and
$0<s<1$.
\par Let us start with the case $s=1$. Then
$$T_2(a)\,= \,\int_{\D
}|f(z)-f(a)|^2 |g^\prime(z)|^2 (1-\vert \varphi
_a(z)\vert^2)\,dA(z).$$ Making the change of variable $w=\varphi
(z)$ in the last integral, we obtain
$$T_2(a)\,= \,\int_{\D }\vert (f\circ \varphi _a)(w)-f(a)\vert
^2\vert (g\circ \varphi _a)^\prime (w)\vert ^2(1-\vert w\vert
^2)\,dA(w).$$ Since $Q_{1, \log , 1/p^\prime }\subset Q_1=BMOA$,
$g\in BMOA$ and then it follows that, for all $a\in \D$, $g\circ
\varphi_a\in BMOA$ and $\rho_*(g\circ \varphi_a)=\rho_*(g)$. This
gives that all the measures $(1-\vert w\vert ^2)\vert (g\circ
\varphi_a)^\prime (w)\vert ^2\,dA(w)$ ($a\in \mathbb D$) are
Carleson measures with constants controlled by $\Vert g\Vert
_{BMOA}^2$. Then, using the Carleson embedding theorem for $H^2$ and
the fact that $B^p$ is continuously embedded in $BMOA$, it follows
that
$$T_2(a)\,\lesssim \Vert g\Vert _{BMOA}^2\Vert f\circ \varphi
_a\,-\-f(a)\Vert_{H^2}^2\,\lesssim
 \Vert g\Vert _{BMOA}^2\Vert f\Vert_{BMOA}^2 \lesssim
 \Vert g\Vert _{BMOA}^2\Vert f\Vert_{B^p}^2.$$
Putting together this, (\ref{muFsCar}), and (\ref{T1BMOA}), we see
that the measure $\mu_{ _{F}}$ is a Carleson measure. This finishes
the proof of part (iii).
\par\medskip To finish the proof of part (iv) we proceed to
estimate $T_2(a)$ assuming that $2<p<\infty $, $0<s<1$,  and
$1-\frac{2}{p}<s$. Notice that $$T_2(a)\,=\,(1-\vert a\vert
^2)^s\int_{\D }\left \vert \frac{f(z)-f(a)}{(1-\overline
a\,z)^s}\right \vert ^2\vert g^\prime (z)\vert ^2(1-\vert z\vert
^2)^s\,dA(z).$$ Since $0<s<1$, $2/p^\prime
>1$, and the measure $(1-\vert z\vert ^2)^s\vert g^\prime (z)\vert
^2\,dA(z)$ is a $2/p^\prime $-logarithmic $s$-Carleson measure,
using Theorem\,\@\ref{Collect-Results}\,\@(v), it follows that
$$T_2(a)\,\lesssim \,(1-\vert a\vert ^2)^s\left (\vert f(a)-f(0)\vert ^2\,+\,\int_{\D }
\left \vert \left (\frac{f(z)-f(a)}{(1-\overline a\,z)^s}\right
)^\prime\right \vert^2(1-\vert z\vert ^2)^s\,dA(z)\right ).$$ The
growth estimate (\ref{grBp}) and simple computations yield
\begin{align*}T_2(a)\,\lesssim \,\, &\Vert f\Vert _{B^p}^2(1-\vert
a\vert ^2)^s\left (\log \frac{2}{1-\vert a\vert ^2}\right
)^{2/p^\prime }\,+\,\int_{\D }\vert f^\prime (z)\vert
^2(1-\vert\varphi _a(z)\vert ^2)^s\,dA(z) \\
\,+\,\,& \int_{\D}\frac{\vert f(z)-f(a)\vert ^2}{\vert 1-\overline
a\,z\vert ^2}(1-\vert \varphi _a(z)\vert ^2)^s\,dA(z) \\
\,\lesssim\,\,&\Vert f\Vert _{B^p}^2\,+\,\int_{\D }\vert f^\prime
(z)\vert ^2(1-\vert\varphi _a(z)\vert
^2)^s\,dA(z)\,+\,\int_{\D}\frac{\vert f(z)-f(a)\vert ^2}{\vert
1-\overline a\,z\vert ^2}(1-\vert \varphi _a(z)\vert ^2)^s\,dA(z).
\end{align*}
By Theorem\,\@\ref{Collect-Results}\,\@(iv), our assumptions on $s$
and $p$ imply that $B^p$ is continuously embedded in $Q_s$. Hence,
$f\in Q_s$. This implies that
$$\int_{\D }\vert f^\prime
(z)\vert ^2(1-\vert\varphi _a(z)\vert ^2)^s\,dA(z)\,\le \, \Vert
f\Vert _{Q_s}^2\lesssim \,\Vert f\Vert _{B^p}^2$$ and that
$$\int_{\D}\frac{\vert f(z)-f(a)\vert ^2}{\vert
1-\overline a\,z\vert ^2}(1-\vert \varphi _a(z)\vert
^2)^s\,dA(z)\,\lesssim \,\Vert f\Vert _{Q_s}^2\lesssim \,\Vert
f\Vert _{B^p}^2,$$ by a result proved by
 Pau and Pel\'{a}ez in
\cite[pp.\,\@551--552]{Pau-Pe-MZ-2009}. Consequently, we have proved
that $T_2(a)\lesssim \,\Vert f\Vert _{B^p}^2$. This, together with
(\ref{muFsCar}) and (\ref{T1BMOA}), shows that $\mu _F$ is an
$s$-Carleson measure as desired. Thus the proof is also finished in
this case.
\end{Pf}
\par\medskip
The case when $1<p\le 2$ and $0<s<1$ remains open. This is so
because if we set $\alpha =2/p^\prime $, then $\alpha \le 1$ and,
hence, $\alpha $ is not in the conditions of
Theorem\,\@\ref{Collect-Results}\,\@(v). We can prove the following
result.
\begin{theorem}\label{Tg1p20s1} Suppose that $1<p\le 2$ and $0<s<1$,
and let $g\in \hol (\D )$. The following statements hold.
\begin{itemize}\item[(i)] If $T_g$ maps $B^p$ into $Q_s$ then $g\in
Q_{s, \log , 1/p^\prime }$.
\item[(ii)] If $\alpha >1/2$ and  $g\in
Q_{s, \log , \alpha }$ then $T_g$ maps $B^p$ into $Q_s$.
\end{itemize}
\end{theorem}
\begin{pf}
(i) follows from part (ii) of Theorem\,\@\ref{TgBpQs}. \par Let us
turn to prove (ii). Suppose that $0<s<1$, $\alpha >1/2$, and $g\in
Q_{s, \log , \alpha }$. Set $$K(g)= \sup _{a\in \D }\left
(\log\frac{2}{1-\vert a\vert }\right )^{2\alpha }\int_{\D}\vert
g^\prime (z)\vert ^2(1-\vert\varphi _a(z)\vert ^2)^s\,dA(z),$$ and
take $f\in B^p$. Set $F=T_g(f)$, we have to prove the $F\in Q_s$ or,
equivalently, that the measure $\mu_{ _{F}}$ defined by
$$d\mu_{ _{F}}(z)\,=\,(1-|z|^2)^s\vert F^\prime(z)|^2\,dA(z)$$ is an $s$-Carleson
measure. Now we argue as in the proof of
Theorem\,\@\ref{TgBpQs}\,\@(iv). For $a\in \D$, we obtain

\begin{equation}\label{muFsCarbis}
\frac{1}{(1-\vert a\vert ^2)^s}\int_{S(a)}
|F^\prime(z)|^2(1-|z|^2)^s\,dA(z)\,\lesssim \,\,
2T_1(a)\,+\,2T_2(a),
\end{equation}
where $T_1(a)$ and $T_2(a)$ are defined as in the proof of
Theorem\,\@\ref{TgBpQs}. Using (\ref{grBp}) and the fact that
$\frac{1}{p^\prime }\le \frac{1}{2}<\alpha $,  we obtain

\begin{equation*}\label{grBpbis}\vert f(a)-f(0)\vert
\lesssim \Vert f\Vert _{B^p}\left (\log \frac{2}{1-\vert a\vert
^2}\right )^{\alpha }.\end{equation*} This yields
\begin{align}\label{T1BMOAbis}T_1(a)
 \lesssim & \,\, \Vert f\Vert
_{B^p}^2\left (\log \frac{2}{1-\vert a\vert^2}\right
 )^{2\alpha }\int_{\D }\vert g^\prime (z)\vert ^2(1-\vert
\varphi _a(z)\vert ^2)^s\,dA(z)  \lesssim \,\,  K(g)\Vert f\Vert
_{B^p}^2.
\end{align}

To estimate  $T_2(a)$, observe that the measure $(1-\vert z\vert
^2)^s\vert g^\prime (z)\vert ^2\,dA(z)$ is a $2\alpha $-logarithmic
$s$-Carleson measure. Since $2\alpha >1$, using Lemma\,\@1 of
\cite{Pau-Pe-MZ-2009}, this implies that the measure $(1-\vert
z\vert ^2)^s\vert g^\prime (z)\vert ^2\,dA(z)$ is a Carleson measure
for $\mathcal D_s$. Then, arguing as in the proof of
Theorem\,\@\ref{TgBpQs}\,\@(iv), we obtain $T_2(a)\lesssim \Vert
f\Vert _{B^p}^2$. This, together with (\ref{T1BMOAbis}) and
(\ref{muFsCarbis}), implies that the measure $\mu_{ _{F}}$ is an
$s$-Carleson measure.
\end{pf}

\par\medskip Regarding the operators $I_g$ and $M_g$ we have the
following results. \begin{theorem}\label{IgMgBpQs} Let $g\in \hol
(\D )$, then:
\begin{itemize}
\item[(1)]
If\, $1<p\le 2$\, and \,$0<s\le 1$\, then: \begin{itemize}
\item[(1a)]\,\, $I_g$ maps $B^p$ into $Q_s$ if and only if $g\in H^\infty $.
\item[(1b)]\,\, If $M_g$\, maps $B^p$ into $Q_s$ then $g\in Q_{s, \log ,
1/p^\prime }\cap H^\infty $. \item[(1c)]\,\, If $g\in Q_{s, \log ,
\alpha }\cap H^\infty $ for some $\alpha
>1/2$ then $M_g$\, maps $B^p$ into $Q_s$.
\end{itemize}
\item[(2)] If\, $2<p<\infty $\, and \,$1-\frac{2}{p}<s\le 1$\, then:
\begin{itemize}
 \item[(2a)]\,\, $I_g$ maps $B^p$ into $Q_s$ if and only if $g\in H^\infty $.
\item[(2b)]\,\, $M_g$ maps $B^p$ into $Q_s$ if and only if $g\in Q_{s, \log ,
1/p^\prime }\cap H^\infty $.\end{itemize}
\item[(3)] If\, $2<p<\infty $\, and \,$0<s\le 1-\frac{2}{p}$\, then:
\begin{itemize} \item[(3a)]\,\, $I_g$ maps $B^p$ into $Q_s$ if and only if $g\equiv 0$.
\item[(3b)]\,\, $M_g$ maps $B^p$ into $Q_s$ if and only if $g\equiv
0$.\end{itemize}
\end{itemize}
\end{theorem}

\begin{Pf}{\it Parts (1) and (2) of Theorem\,\@\ref{IgMgBpQs}.}
 Using
Proposition\,\@\ref{MXY-Conf} it follows that if either $I_g$ or
$M_g$ maps $B^p$ into $Q_s$ for any pair $(s, p)$ with $0<s<\infty $
and $1<p<\infty $ then $g\in H^\infty $.
\par Suppose now that $s$ and $p$ are in the conditions of (1) or
(2) and that $g\in H^\infty $. Take $f\in B^p$. We have to prove
$I_g(f)\in Q_s$ or, equivalently, that the measure
\begin{equation}\label{Igsc}(1-\vert z\vert ^2)^s\vert f^\prime (z)\vert ^2\vert
g(z)\vert ^2\,dA(z) \text{is an $s$-Carleson measure.}\end{equation}
Using (\ref{BpQsple2}) and (\ref{BpQsp>2}), we see that $B^p\subset
Q_s$. Hence $f\in Q_s$ which is the same as saying that $(1-\vert
z\vert ^2)^s\vert f^\prime (z)\vert ^2\,dA(z)$ is an $s$-Carleson
measure. This and the fact that $g\in H^\infty $ trivially yield
(\ref{Igsc}). Thus (1a) and (2a) are proved. Then (1b), (1c), and
(2b) follow using Proposition\,\@\ref{MXY-Conf}, the fact that if
two of the operators $T_g$, $I_g$, $M_g$ map $B^p$ into $Q_s$ so
does the third one, Theorem\,\@\ref{TgBpQs}, and
Theorem\,\@\ref{Tg1p20s1}.
\end{Pf}
\par\medskip
In order to prove Theorem\,\@\ref{IgMgBpQs}\,\@(3), for $2<p<\infty
$ we shall consider the function $F_p$ defined by
\begin{equation}\label{Fp}F_p(z)\,=\,\sum_{k=1}^\infty
\frac{1}{k^{1/2}2^{k/p}}\,z^{2^k},\quad z\in \D .\end{equation}
Using \cite[Corollary\,\@7]{AC} or \cite[Theorem\,\@6]{AXZ}, we see
that $F_p\in B^p$ and $F_p\notin Q_{1-\frac{2}{p}}$. Hence
\begin{equation}\label{Fpnotin}F_p\in B^p\setminus Q_s,\quad
0<s\le 1-\frac{2}{p}, \,\,\,\,2<p<\infty .\end{equation}
\par Let us estimate the integral means $M_2(r, F_p^\prime )$. We
have
$$zF_p^\prime (z)\,=\,\sum_{k=1}^\infty 2^{k/p^\prime
}k^{-1/2}z^{2^k},\quad z\in \D$$ and, hence,
$$M_2(r, F_p^\prime )^2\,\gtrsim \,\sum_{k=1}^\infty 2^{2k/p^\prime
}k^{-1}\,r^{2^{k+1}},\quad 0<r<1.$$ Set $r_n=1-2^{-n}$ ($n=1, 2,
\dots $). Then \begin{align*}&M_2(r_n, F_p^\prime )^2\,\gtrsim
\,\sum_{k=1}^\infty 2^{2k/p^\prime }k^{-1}\,r_n^{2^{k+1}}\\
&\,\gtrsim \, 2^{2n/p^\prime }n^{-1}\,r_n^{2^{n+1}}\,\gtrsim
\,2^{2n/p^\prime }n^{-1}\,\asymp\,\frac{1}{(1-r_n)^{2/p^\prime }\log
\frac{2}{1-r_n}},\quad n=1, 2, \dots .\end{align*} This readily
yields
\begin{equation}\label{M2Fp}M_2(r, F_p^\prime )^2\,\gtrsim \frac{1}{(1-r)^{2/p^\prime }\log
\frac{2}{1-r}},\quad 0<r<1.\end{equation}
\par\medskip
\begin{Pf}{\it part (3) of Theorem\,\@\ref{IgMgBpQs}.} Suppose that
$2<p<\infty $\, and \,$0<s\le 1-\frac{2}{p}$ and  $g\in \hol (\D )$
is not identically zero.
\par Suppose first that either $I_g$ or $M_g$ maps $B^p$ into $Q_s$. We know
that then $g\in H^\infty $ and then, by Fatou's theorem and the
Riesz uniqueness theorem, we know  that $g$ has a finite
non-tangential limit $g(e^{i\theta })$ for almost every $\theta \in
[0, 2\pi ]$ and that $g(e^{i\theta })\neq 0$ for almost every
$\theta $. Then it follows that there exist $C>0$, $r_0\in (0, 1)$,
and a measurable set $E\subset [0, 2\pi ]$ whose Lebesgue measure
$\vert E\vert $ is positive such that
\begin{equation}\label{gE}\vert g(re^{i\theta })\vert \ge C,\quad
\theta \in E,\,\,\,\,r_0<r<1.\end{equation} Since $F_p$ is given by
a power series with Hadamard gaps, Lemma\,\@6.\,\@5 in
\cite[Vol.\,\@1, p.\,\@203]{Zy} implies that
\begin{equation}\label{FpallE}\int_E\,\vert F_p^\prime (re^{i\theta
})\vert ^2\,d\theta \,\asymp \,M_2(r, F_p^\prime )^2,\quad
0<r<1.\end{equation}
\par\medskip
Using the fact that $s\le 1-\frac{2}{p}$, (\ref{gE}),
(\ref{FpallE}), and (\ref{M2Fp}), we obtain
\begin{align}\label{intFpg}
&\int_0^1\,(1-r)^sM_2(r, F_p^\prime\,g)^2\,dr\,\ge
\,\int_{r_0}^1\,(1-r)^{1-\frac{2}{p}}M_2(r,
F_p^\prime\,g)^2\,dr\nonumber \\ \gtrsim &
\int_{r_0}^1\,(1-r)^{1-\frac{2}{p}}\int _E\vert F_p^\prime
(re^{i\theta })\vert ^2\,\vert g(re^{i\theta })\vert
^2\,d\theta\,dr\, \gtrsim\,\int_{r_0}^1\,(1-r)^{1-\frac{2}{p}}\int
_E\vert F_p^\prime (re^{i\theta })\vert ^2\,d\theta\,dr\nonumber \\
\gtrsim & \int_{r_0}^1\,(1-r)^{1-\frac{2}{p}}\,M_2(r, F_p^\prime
)^2\,dr\,\gtrsim \int_{r_0}^1\frac{dr}{(1-r)\log \frac{2}{1-r}} =
\infty .
\end{align}

\par\medskip
If we assume that $I_g$ maps $B^p$ into $Q_s$ then $I_g(F_p)\in Q_s$
and then, using \cite[Proposition\,\@3.\,\@1]{AGW}, it follows that
$$\int_0^1\,(1-r)^sM_2(r, F_p^\prime\,g)^2\,dr\,<\infty .$$ This is
in contradiction with (\ref{intFpg}).
\par Assume now that $M_g$ maps $B^p$ into $Q_s$.
Since $1$ and $F_p$ belong to $B^p$, we have that $g$ and $F_p\,g$
belong to $Q_s$ and then, by \cite[Proposition\,\@3.\,\@1]{AGW},
\begin{equation}\label{int-M2g}\int_0^1(1-r)^s\,M_2(r, g^\prime
)^2\,dr\,<\,\infty \end{equation} and
\begin{equation}\label{int-prod} \int_0^1(1-r)^s\,M_2(r,
\,(F_pg)^\prime )^2\,dr\,<\,\infty .
\end{equation} Notice that
$F_p\in H^\infty $ and then
\begin{equation*}M_2(r, F_p\,g^\prime
)\,\lesssim \,M_2(r, g^\prime ),\quad 0<r<1.\end{equation*} This and
(\ref{int-M2g}) imply that
\begin{equation}\label{int-Fpgprime}\int_0^1(1-r)^s\,M_2(r,
F_p^\prime\,g)^2\,dr\,<\,\infty .\end{equation}  We have arrived to
a contradiction because it is clear that (\ref{intFpg}) and
(\ref{int-Fpgprime}) cannot be simultaneously true.
\end{Pf}
\par\medskip In the other direction we have the following result.
\begin{theorem}\label{QsBP} Suppose that $0<s<\infty $ and $1\le
p<\infty $ and let $g\in \hol (\D )$. Then the following conditions
are equivalent
\begin{itemize} \item[(i)] $M_g$ maps $Q_s$ into $B^p$.
\item[(ii)] $g\equiv 0$.
\end{itemize}
\end{theorem}
\begin{pf} Suppose that $g\not\equiv 0$. Choose an increasing sequence $\{
r_n\} _{n=1}^\infty \subset (0, 1)$ with $\,\lim \{r_n\}=1$\, and a
sequence $\{ \theta _n\} _{n=1}^\infty \subset [0, 2\pi ]$ such that
$$\vert g(r_ne^{i\theta _n})\vert =M_\infty (r_n, g),\quad n=1, 2,
\dots .$$ For each $n$ set $$f_n(z)\,=\,\log \frac{1}{1-e^{-i\theta
_n}z},\quad z\in \D .$$ Notice that $M(r_1, g)>0$ and that the
sequence $\{ M(r_n, g)\} $ is increasing. Set $$f_n(z)=\log
\frac{1}{1-e^{-i\theta _nz}},\quad z\in \D,\quad n=1, 2, \dots .$$
We have that $f_n\in Q_{s}$ for all $n$ and
$$\Vert f_n\Vert _{Q_s}\asymp
1.$$
\par Assume that $M_g$ maps $Q_s$ into $B^p$. Then, by the closed graph theorem, $M_g$ is
bounded operator from $Q_s$ into $B^p$. Hence the sequence $\{
g\,f_n \} _{n=1}^\infty $ is a bounded sequence on $B^p$, that is,
$$\Vert g\,f_n\Vert _{B^p}\lesssim 1.$$ Then it follows that
\begin{align*}M(r_1, g)\log \frac{1}{1-r_n}\le& M(r_n, g)\log
\frac{1}{1-r_n}=\vert g(r_ne^{i\theta _n})f_n(r_ne^{i\theta
_n})\vert\\ \lesssim & \Vert g\,f_n\Vert_{B^p}\left (\log
\frac{1}{1-r_n}\right )^{1/p^\prime }\lesssim \left (\log
\frac{1}{1-r_n}\right )^{1/p^\prime }.\end{align*} This is a
contradiction.
\end{pf}

\section{Multipliers and integration operators between $Q_s$ spaces}
As we mentioned above the space of multipliers $M(\mathcal
B)=M(Q_s)$ ($s>1$) was characterized by Brown and Shields in
\cite{BS}. Ortega and F\`{a}brega \cite{OF} characterized the space
$M(BMOA)=M(Q_1)$. Pau and Pel\'{a}ez \cite{Pau-Pe-MZ-2009} and Xiao
\cite{X-2008} characterized the spaces $M(Q_s)$ ($0<s<1$) closing a
conjecture formulated in \cite{X-Pac}. Indeed, Theorem\,\@1 of
\cite{Pau-Pe-MZ-2009} and Theorem\,\@1.\,\@2 of \cite{X-2008} assert
the following.
\begin{other}\label{mult-Qs}
Suppose that $0<s\le 1$ and let $g$ be an analytic function in the
unit disc $\D $. Then:
\begin{itemize}
\item[(i)] $T_g$ maps $Q_s$ into itself if and only if $g\in
Q_{s, \log , 1}$.
\item[(ii)] $I_g$ maps $Q_s$ into itself if and only if $g\in
H^\infty $.
\item[(ii)] $M_g$ maps $Q_s$ into itself if and only if $g\in Q_{s, \log ,
1}\cap H^\infty $.
\end{itemize}
\end{other}
\par\medskip We shall prove the following results.
\begin{theorem}\label{Qs1les2} Suppose that $0<s_1\le s_2\le 1$ and
let $g\in \hol (\D )$. Then:
\begin{itemize}
\item[(i)] $T_g$ maps $Q_{s_1}$ into $Q_{s_2}$ if and only if $g\in
Q_{{s_2},\log ,1}$.
\item[(ii)] $I_g$ maps $Q_{s_1}$ into $Q_{s_2}$ if and only if $g\in
H^\infty $.
\item[(iii)] $M_g$ maps $Q_{s_1}$ into $Q_{s_2}$ if and only if $g\in
Q_{{s_2},\log ,1}\cap H^\infty $.
\end{itemize}
\end{theorem}
\par\medskip
\begin{theorem}\label{Qs1mas2} Suppose that $0<s_1<s_2\le 1$ and
let $g\in \hol (\D )$. Then the following conditions are equivalent:
\begin{itemize}
\item[(i)] $I_g$ maps $Q_{s_2}$ into $Q_{s_1}$.
\item[(ii)] $M_g$ maps $Q_{s_2}$ into $Q_{s_1}$.
\item[(iii)] $g\equiv 0$.
\end{itemize}
\end{theorem}

\begin{Pf}{\it Theorem\,\@\ref{Qs1les2}.} For $a\in \mathbb D$ we
set
$$h_a(z)\,=\,\log\frac{2}{1-\overline a\,z},\quad z\in \D.$$
Then $h_a\in Q_{s_1}$ for all $a\in \D$ and
\begin{equation}\label{normha}\Vert h_a\Vert _{Q_{s_1}}\asymp
1.\end{equation}
\par $\bullet $ If $T_g$ maps $Q_{s_1}$ into $Q_{s_2}$ then $T_g$ is
a bounded operator from $Q_{s_1}$ into $Q_{s_2}$. Using this and
(\ref{normha}), it follows that for all $a\in \D$ the measure
$(1-\vert z\vert^2)^{s_2}\vert g^\prime (z)\vert ^2\vert h_a(z)\vert
^2\,dA(z)$ is an $s_2$-Carleson measure and that
\begin{equation}\label{s2car}\int_{S(a)}(1-\vert z\vert
^2)^{s_2}\vert g^\prime (z)\vert ^2\vert h_a(z)\vert
^2\,dA(z)\,\lesssim (1-\vert a\vert ^2)^{s_2},\quad a\in \D
.\end{equation} Since $$\vert h_a(z)\vert \asymp \log
\frac{2}{1-\vert a\vert ^2},\quad z\in S(a),$$ (\ref{s2car}) implies
that
$$\left (\log\frac{2}{1-\vert a\vert ^2}\right
)^2\int_{S(a)}(1-\vert z\vert^2)^{s_2}\vert g^\prime (z)\vert
^2\,dA(z)\,\lesssim (1-\vert a\vert ^2)^{s_2}.$$ This is the same as
saying that the measure $(1-\vert z\vert ^2)^{s_2}\vert g^\prime
(z)\vert ^2\,dA(z)$ is a $2$-logarithmic $s_2$-Carleson measure or,
equivalently, that $g\in Q_{s_2, \log ,1}$.
\par If $g\in Q_{s_2, \log ,1}$ then, by Theorem\,\@\ref{mult-Qs},
$T_g$ maps $Q_{s_2}$ into itself. Since $Q_{s_1}\subset Q_{s_2}$, it
follows trivially that $T_g$ maps $Q_{s_1}$ into $Q_{s_2}$. Hence
(i) is proved
\par $\bullet $ Proposition\,\@\ref{MXY-Conf} shows that if $I_g$
maps $Q_{s_1}$ into $Q_{s_2}$ then $g\in H^\infty $.
\par Conversely, suppose that $g\in H^\infty $. In order to prove
that $I_g$ maps $Q_{s_1}$ into $Q_{s_2}$, we have to prove that for
any $f\in Q_{s_1}$ the measure $(1-\vert z\vert ^2)^{s_2}\vert
g(z)\vert ^2\vert f^\prime (z)\vert ^2\,dA(z)$ is an $s_2$-Carleson
measure. So, take $f\in Q_{s_1}$. Then $(1-\vert z\vert
^2)^{s_1}\vert f^\prime (z)\vert ^2\,dA(z)$ is an $s_1$-Carleson
measure. Then it follows that \begin{align*}&\int_{S(a)}(1-\vert
z\vert ^2)^{s_2}\vert g(z)\vert ^2\vert f^\prime (z)\vert
^2\,dA(z)\\ \le &\,\Vert g\Vert_{H^\infty }^2(1-\vert a\vert
^2)^{s_2-s_1}\int_{S(a)}(1-\vert z\vert ^2)^{s_1}\vert f^\prime
(z)\vert ^2\,dA(z)\\ \lesssim &\, (1-\vert a\vert
^2)^{s_2}.\end{align*} This shows that $(1-\vert z\vert
^2)^{s_2}\vert g(z)\vert ^2\vert f^\prime (z)\vert ^2\,dA(z)$ is an
$s_2$-Carleson measure as desired, finishing the proof of (ii).
\par $\bullet $ If $M_g$ maps $Q_{s_1}$ into $Q_{s_2}$ then, Proposition\,\@\ref{MXY-Conf}, $g\in
H^\infty $. Then (i) implies that $I_g$ maps $Q_{s_1}$ into
$Q_{s_2}$. Since $M_g(f)\,=I_g(f)\,+\,T_g(f)\,+f(0)g(0)$, it follows
that $T_g$ maps $Q_{s_1}$ into $Q_{s_2}$. Then (i) yields $g\in
Q_{s_2, \log , 1}$. Then we have that  $g\in Q_{s_2, \log , 1}\cap
H^\infty $.
\par Conversely, if $g\in
Q_{s_2, \log , 1}\cap H^\infty $ then (i) and (ii) immediately give
that both $T_g$ and $I_g$ map $Q_{s_1}$ into $Q_{s_2}$ and then so
does $M_g$.
\end{Pf}
\par\medskip Some results from \cite{AGW} will be used to prove
Theorem\,\@\ref{Qs1mas2}. As we have already noticed if $0<s\le 1$
and  $f\in Q_s$ then $\int_0^1(1-r)^sM_2(r, f^\prime
)^2\,dr\,<\infty $. Using ideas from \cite{G-der-boun}, Aulaskari,
Girela and Wulan \cite[Theorem\,\@3.\,\@1]{AGW} proved that this
result is sharp in a very strong sense.
\begin{other}\label{sharp} Suppose that $0<s\le 1$ and let $\varphi
$ be a positive increasing function defined in $(0, 1)$ such that
$$\int_0^1(1-r)^s\,\varphi (r)^2\,dr\,<\infty .$$ Then there exists a
function $f\in Q_s$ given by a power series with Hadamard gaps such
that $M_2(r, f^\prime )\ge \varphi (r)$ for all $r\in (0, 1)$.
\end{other}
\par\medskip
\begin{Pf}{\it Theorem\,\@\ref{Qs1mas2}.}
Suppose that $g\not\equiv 0$ and that either $I_g$ or $M_g$ maps
$Q_{s_2}$ into $Q_{s_1}$. By Proposition\,\@\ref{MXY-Conf}, $g\in
H^\infty $ and then it follows that there exist $C>0$, $r_0\in (0,
1)$, and a measurable set $E\subset [0, 2\pi ]$ whose Lebesgue
measure $\vert E\vert $ is positive such that $$\vert g(re^{i\theta
})\vert \ge C,\quad \theta \in E,\quad r_0<r<1.$$ \par $\bullet $
Suppose that $I_g$ maps $Q_{s_2}$ into $Q_{s_1}$. Then we use
Theorem\,\@\ref{sharp} to pick a function $F\in Q_{s_2}$ given by a
power series with Hadamard gaps so that
\begin{equation}\label{M2rFrime-big}M_2(r, F^\prime )\,\ge
\,\frac{1}{(1-r)^{(1+s_1)/2}},\quad 0<r<1.\end{equation} Since
$I_g(F)\in Q_{s_1}$,
\begin{equation}\label{IntF}\int_0^1(1-r)^{s_1}M_2(r, F^\prime
g)^2\,dr \,<\,\infty .\end{equation} However, using Lemma\,\@6.\,\@5
in \cite[Vol.\,\@1, p.\,\@203]{Zy} and (\ref{M2rFrime-big}), it
follows that
\begin{align*}\int_0^1(1-r)^{s_1}M_2(r, F^\prime
g)^2\,dr &\, \gtrsim \,\int_{r_0}^1(1-r)^{s_1}\int_E\vert F^\prime
(re^{i\theta })\vert ^2\vert g(re^{i\theta })\vert ^2\,d\theta \,dr
\\ & \,\gtrsim \,\int_{r_0}^1(1-r)^{s_1}\int_E\vert F^\prime
(re^{i\theta })\vert^2\,d\theta \,dr \\ & \asymp
\,\int_{r_0}^1(1-r)^{s_1}M_2(r, F^\prime )^2\,dr \\ & \,\gtrsim
\,\int_{r_0}^1(1-r)^{-1}\,dr\\ & \,= \,\infty .
\end{align*}
This is in contradiction with (\ref{IntF}).
\par \par $\bullet $
Suppose now that $M_g$ maps $Q_{s_2}$ into $Q_{s_1}$. Take
$\varepsilon >0$ with $s_2-s_1-\varepsilon >0$ and use
Theorem\,\@\ref{sharp} to pick a function $H\in Q_{s_2}$ given by a
power series with Hadamard gaps so that
\begin{equation}\label{M2rHrime-big}M_2(r, H^\prime )\,\ge
\,\frac{1}{(1-r)^{(1+s_1+\varepsilon )/2}},\quad
0<r<1.\end{equation} Since $gH\in \,Q_{s_1}$ we have that
\begin{equation}\label{prod-prime}\int_0^1(1-r)^{s_1}M_2(r, g^\prime
H\,+\,gH^\prime )^2\,dr\,<\,\infty .\end{equation} Using
Lemma\,\@6.\,\@5 in \cite[Vol.\,\@1, p.\,\@203]{Zy} and
(\ref{M2rHrime-big}), we obtain as above that
\begin{align}\label{Hprimeg}\int_0^1(1-r)^{s_1+\varepsilon }M_2(r,
H^\prime g)^2\,dr\,\gtrsim \, & \int_{r_0}^1(1-r)^{s_1+\varepsilon
}\int_E\,\vert H^\prime (re^{i\theta })\vert^2\,d\theta \,dr \nonumber \\
\,\gtrsim \, & \int_{r_0}^1(1-r)^{s_1+\varepsilon }M_2(r, H^\prime
)^2\,dr \nonumber  \\
\,\gtrsim \, & \int_{r_0}^1\frac{dr}{1-r} \nonumber \\ \,= \, &
\infty .\end{align} Notice that $g\in Q_{s_1}$. Using this and the
fact that $$\vert H(z)\vert \lesssim \log \frac{2}{1-\vert z\vert
},\quad z\in \D,$$ it follows that
\begin{align}\label{Hgprimemmm}&\int_0^1(1-r)^{s_1+\varepsilon }M_2(r, Hg^\prime
)^2\,dr\,\lesssim \int_0^1(1-r)^{s_1+\varepsilon }\left (\log
\frac{2}{1-r}\right )^2M_2(r, g^\prime )^2\,dr \nonumber \\ \, &
\lesssim \int_0^1(1-r)^{s_1+\frac{\varepsilon }{2}}M_2(r, g^\prime
)\,dr\,<\,\infty .\end{align} We have arrived to a contradiction
because (\ref{prod-prime}), (\ref{Hprimeg}), and (\ref{Hgprimemmm})
cannot hold simultaneously.
\end{Pf}

\par\medskip
\begin{remark} The implication (ii)\, $\Rightarrow$\, (iii) in
Theorem\,\@\ref{Qs1mas2} was obtained by Pau and Pel\'{a}ez
\cite[Corollary\,\@4]{Pau-Pe-2011} using the fact that there exists
a function $f\in Q_{s_2}$, $f\not\equiv 0$, whose sequence of zeros
is not a $Q_{s_1}$-zero set.
\par This idea gives also the following:
\par $$M(\mathcal B, Q_s)\,=\,\{ 0\},\quad 0<s\le 1.$$
\par\medskip Indeed, it is well known that there exists a function
$f\in \mathcal B$, $f\not\equiv 0$, whose sequence of zeros does not
satisfy the Blaschke condition \cite{ACP,GNW}. If $g\not\equiv 0$
were a multiplier from $\mathcal B$ into $Q_s$ for some $s\le 1$
then the sequence of zeros of $fg$ would satisfy the Blaschke
condition. But this is not true because all the zeros of $f$ are
zeros of\, $gf$.
\end{remark}

\par\medskip
\section{Some further results}
\par The inner-outer factorization of functions in the Hardy spaces
plays an outstanding role in lots of questions. In many cases the
outer factor $O_f$ of $f$ inherits properties of $f$. Working in
this setting the following concepts arise as natural and quite
interesting.
\par A subspace $X$ of $H^1$ is said to have the
$f$-property (also called the property of division by inner
functions) if $h/I\in X$ whenever $h\in X$ and $I$ is an inner
function with $h/I\in H^1$.
\par Given $v\in L^\infty (\partial \D )$, the Toeplitz operator
$T_v$ associated with the symbol $v$ is defined by
$$T_vf(z)\,=\,P(vf)(z)\,=\,\frac{1}{2\pi i}\int_{\partial
\D}\frac{v(\xi )f(\xi )}{\xi -z}\,d\xi,\quad f\in H^1,\quad z\in \D
.$$ Here, $P$ is the Szeg\"{o} projection.
\par A subspace $X$ of $H^1$ is said to have the $K$-property if
$T_{\overline \psi }(X)\subset X$ for any $\psi \in H^\infty $.
\par\medskip These notions were introduced by Havin in \cite{Havin}.
It was also pointed out in \cite{Havin} that the $K$-property
implies the $f$-property: indeed, if $h\in H^1$, $I$ is inner and
$h/I \in H^1$ then $h/I=T_{\overline I} h$.
\par In addition to the Hardy spaces $H^p$ ($1<p<\infty $) many other
spaces such as the Dirichlet space \cite{Havin,kor:outer}, several
spaces of Dirichlet type including all the Besov spaces $B^p$
($1<p<\infty $) \cite{dya:h-s,Dya:hpw,Dya:gar,kor-faivy73}, the
spaces $BMOA$ and $VMOA$ \cite{he:kprop}, and the $Q_s$ spaces
($0<s<1$) \cite{DyG} have the $K$-property. The Hardy space $H^1$,
$H^\infty $ and $VMOA\cap H^\infty $ are examples of spaces which
have the $f$-property bur fail to have the $K$-property
\cite{he:kprop}. \par The first example of a subspace of $H^1$ not
possessing the $f$-property is due to Gurarii \cite{Gu} who proved
that the space of analytic functions in $\D $ whose sequence of
Taylor coefficients is in $\ell ^1$ does not have the $f$-property.
Anderson \cite{an:di} proved that the space $\mathcal B_0\cap
H^\infty $ does not have the $f$-property. Later on it was proved in
\cite{GGP} that if $1\le p<\infty $ then $H^p\cap \mathcal B$ fails
to have the $f$-property also.
\par\medskip Since as we have already mentioned the Besov spaces
$B^p$ ($1<p<\infty $) and the $Q_s$ spaces ($0<s\le 1$) have the
$K$-property (and, also, the $f$-property), it seems natural to
investigate whether the spaces of multipliers and the spaces $Q_{s,
\log , \alpha }$ that we have considered in our work have also these
properties. We shall prove the following results.
\begin{theorem}\label{mult-f-prop} The spaces of multipliers $M(B^p,
Q_s)$ \,$(0<s\le 1, 1\le p<\infty )$, $M(Q_{s_1}, Q_{s_2})$
\,$(0<s_1, s_2\le 1)$, and $M(B^p, B^q)$\, $(1\le p, q<\infty )$
have the $f$-property.
\end{theorem}
\begin{theorem}\label{QalphalogK} For $\alpha >0$ and $0<s<1$ the
space $Q_{s, \log ,\alpha }$ has the $K$-property.
\end{theorem}
\par\medskip Theorem\,\@\ref{mult-f-prop} follows readily from the
following result.
\begin{lemma}\label{f-pro-mult} Let $X$ and $Y$ be to Banach spaces
of analytic functions which are continuously contained in $H^1$.
Suppose that $X$ contains the constants functions and that $Y$ has
the $f$-property. Then the space of multipliers $M(X, Y)$ also has
the $f$-property.
\end{lemma}

\begin{pf} Since $X$ contains the constants functions $M(X, Y)\subset Y\subset H^1$. \par Suppose that $F\in M(X, Y)$, $I$ is an
inner function, and $F/I\in H^1$. Take $f\in X$. Then $fF\in
Y\subset H^1$ and then $fF/I\in H^1$. Since $Y$ has the
$f$-property, it follows that $fF/I\in Y$. Thus, we have proved that
$F/I\in M(X, Y)$. \end{pf}

\par\medskip Theorem\,\@\ref{QalphalogK} will follows from a
characterization of the spaces $Q_{s, \log ,\alpha }$ in terms of
pseudoanalytic continuation. We refer to Dyn'kin's paper
\cite{Dyn:pse} for similar descriptions of classical smoothness
spaces, as well as for other important applications of the
pseudoanalytic extension method.

\par
Let, $\D_-$ denotes the region $\{ z\in \mathbb C : \vert z\vert
>1\}$, and write
$$z^* \ig 1/{\overline z}, \quad z\in \mathbb C \setminus \{0\}.$$
We shall need the Cauchy-Riemann operator
$$\dbar\, =\,\frac{\partial }{\partial \overline z}\,\ig \,\frac{1}{2}\left (
\frac{\partial }{\partial x}+i\frac{\partial }{\partial y}\right
),\quad z=x+iy.$$ The following result is an extension of
\cite[Theorem\,\@1]{DyG}.
\begin{theorem}\label{pseudo}
Suppose that $0<s<1$, $\alpha>0$, and  $f\in
\displaystyle{\cap_{0<q<\infty }H\sp q}$. Then the following
conditions are equivalent. \begin{itemize} \item[(i)] $\quad f\in
Q_{s,\log,\alpha}$.
\item
[(ii)] $\quad \displaystyle\sup_{\vert a\vert <1}
\left(\log\frac{2}{1-|a|}\right)^{2\alpha}\int_\D \vert
f'(z)\vert\sp 2\left (\frac{1}{\vert \varphi _a(z)\vert \sp
2}-1\right )\sp s\, dA(z)<\infty .$
\item
[(iii)] There exists a function $F\in C\sp 1(\D _-)$ satisfying
\begin{align*}&
F(z)=O(1),\quad \hbox{as $z\to\infty $},\\
&\lim_{r\to 1\sp +}F(re\sp {i\theta })= f(e\sp {i\theta }),\quad
\hbox{a.e. and in $L\sp {q}([-\pi ,\pi ])$ for all $q\in
[1, \infty )$ },\\
&\sup_{\vert a\vert
<1}\left(\log\frac{2}{1-|a|}\right)^{2\alpha}\int_{\D
_-}\left\vert\overline{\partial} F(z)\right\vert \sp 2\left ( \vert
\varphi _a(z)\vert \sp 2-1\right )\sp s\, dA(z)<\infty .\end{align*}
\end{itemize}
\end{theorem}
\par Theorem\,\@\ref{pseudo} can be proved following the arguments
used in the proof of \cite[Theorem\,\@1]{DyG}, we omit the details.
Once Theorem\,\@\ref{pseudo} is established, noticing that if
$\alpha
>0$ and $0<s<1$ then $Q_{s, \log \alpha }\subset Q_s\subset BMOA$,
Theorem\,\@\ref{QalphalogK} can be proved following the steps in the
proof of \cite[Theorem\,\@2]{DyG}. Again, we omit the details.

\par\medskip {\bf Acknowledgements.} We wish to thank the referees for
reading carefully the paper and making a number of nice suggestions
to improve it.
\medskip



\begin{thebibliography}{100}

%
\bibitem{AlCi} A.~Aleman and J.~A.~Cima,
{\em{An integral operator on $H\sp p$ and Hardy's inequality}}, J.
Anal. Math. \textbf{85} (2001), 157--176.
%
\bibitem{ADMV} A.~Aleman, P.~L.~Duren, M.~J.~Mart\'{\i}n and D.~Vukoti\'c, \emph{Multiplicative isometries and isometric zero-divisors},
Canad. J. Math. \textbf{62} (2010), no. 5, 961–-974.
%
\bibitem{AlSim} A.~Aleman and A. Simbotin,
\emph{Estimates in M\"{o}bius invariant spaces of analytic
functions}, Complex Var. Theory Appl. \textbf{49} (2004), no. 7-9,
487–-510.
%
\bibitem{AlSis-CVTA} A.~Aleman and A.~G.~Siskakis,
\emph{An integral operator on $H^p$}, Complex Variables Theory Appl.
\textbf{28} (1995), no. 2, 149–-158.
%
\bibitem{AlSis-Indiana} A.~Aleman and A.~G.~Siskakis,
\emph{Integration operators on Bergman spaces}, Indiana Univ. Math.
J. \textbf{46} (1997), no. 2, 337--356.
%
\bibitem{an:di}
 J.~M.~Anderson,
\emph{On division by inner factors},
 Comment. Math. Helv. \textbf{54} (1979), no.~2, 309--317.
%
\bibitem{ACP}
J.~M.~Anderson, J.~Clunie and Ch.~Pommerenke, \emph{On Bloch
functions and normal functions}, J. Reine Angew. Math. \textbf{270}
(1974), 12--37.
%
\bibitem{AFP}
J.~Arazy, S.~D.~Fisher and J.~Peetre, \emph{M\"obius invariant
function spaces}, J. Reine Angew. Math. \textbf{363} (1985),
110--145.
%
\bibitem{ARS} N.~Arcozzi, R.~Rochberg and E.~Sawyer, \textit{Carleson measures for
analytic Besov spaces.} Rev. Mat. Iberoamericana\/ \textbf{18}
(2002), no.~2, 443--510.
%
\bibitem{AC} R.~Aulaskari and G.~Csordas, \emph{Besov spaces and the $Q_{q,0}$ classes},
Acta Sci. Math. (Szeged) \textbf{60} (1995), no. 1--2, 31-–48.
%

\bibitem{AGW} R.~Aulaskari, D.~Girela and H.~Wulan, \emph{Taylor coefficients and mean growth of the derivative of $Q_p$ functions},
 J. Math. Anal. Appl. \textbf{258} (2001), no. 2, 415–-428.
 %
\bibitem{au-la94}
R.~Aulaskari and P.~Lappan, \emph{Criteria for an analytic function
to be Bloch and a harmonic or meromorphic function to be normal},
Complex Analysis and its Applications (Harlow), Pitman Research
Notes in Math, vol.~305, Longman Scientific and Technical, 1994,
136--146.
%
\bibitem{ASX} R.~Aulaskari, D.~A.~Stegenga and J.~Xiao, \emph{Some subclasses of BMOA and their characterization in terms of Carleson measures},
 Rocky Mountain J. Math. \textbf{26} (1996), no. 2, 485–-506.
%
\bibitem{AXZ} R.~Aulaskari, J.~Xiao and R.~Zhao, \emph{On subspaces and subsets of BMOA and UBC },
Analysis \textbf{15} (1995), no. 2, 101–-121.
%
\bibitem{BS} L.~Brown and A.~L.~Shields, \emph{Multipliers and cyclic vectors in the Bloch
space}, Michigan Math. J. \textbf{38} (1991), no.~1, 141--146.
%
\bibitem{DGV1}
J.~J.~Donaire, D.~Girela and D.~Vukoti\'c, \emph{On univalent
functions in some M\"obius invariant spaces}, J. Reine Angew. Math.
\textbf{553} (2002), 43--72.
%
\bibitem{DGV2}
J.~J.~Donaire, D.~Girela and D.~Vukoti\'c, \emph{On the growth and
range of functions in M\"{o}bius invariant spaces},  J. Anal. Math.
\textbf{112} (2010), 237–-260.
%

\bibitem{D}
P.~L.~Duren, \emph{Theory of $H\sp{p}$spaces}, Academic Press, New
York-London, 1970. Reprint: Dover, Mineola-New York, 2000.
%
\bibitem{DS}
P.~L.~Duren \and A.~P.~Schuster, \emph{Bergman Spaces\/}, Math.
Surveys and Monographs, Vol. 100,  American Mathematical Society,
Providence, Rhode Island, 2004.
%
\bibitem{dya:h-s}
 K.~M.~Dyakonov,
 \emph{Factorization of smooth analytic functions via
 Hilbert-Schmidt operators},
 (in Russian),
 Algebra i Analiz \textbf{8}(1996), no.~4, 1--42.
 {\sl English translation in St. Petersburg Math.~J.
 \textbf{8} (1997), no.~4, 543--569.}

\bibitem{Dya:hpw}
 K.~M.~Dyakonov,
 \emph{Equivalent norms on Lipschitz-type spaces of holomorphic
 functions},
 Acta. Math. \textbf{178} (1997), 143--167.
\bibitem{Dya:gar}
 K.M. Dyakonov,
\emph{Holomorphic functions and quasiconformal mappings with smooth
 modulii},
 Adv. in Math. \textbf{187} (2004), 146--172.

%
\bibitem{DyG} K.~M.~Dyakonov and D.~Girela, \emph{On $Q_p$ spaces and pseudoanalytic extension},
Ann. Acad. Sci. Fenn. Math. \textbf{25} (2000), no. 2, 477–-486.
%
\bibitem{Dyn:pse}
 E.~M.~Dyn'kin,
 \emph{The pseudoanalytic extension},
 J. Anal. Math. \textbf{60} (1993), 45--70.
%
\bibitem{GaGiMa-CAOT} P.~Galanopoulos, D.~Girela and M.~J.~Mart\'{\i}n, \emph{
Besov spaces, multipliers and univalent functions}, Complex Anal.
Oper. Theory \textbf{7} (2013), no. 4, 1081–-1116.
%
\bibitem{GaGiPe-TAMS-2011} P.~Galanopoulos, D.~Girela and J.~A.~Pel\'{a}ez, \emph{Multipliers and integration operators on Dirichlet spaces},
Trans. Amer. Math. Soc. \textbf{363} (2011), no. 4, 1855–-1886.
%
\bibitem{G-der-boun} D.~Girela, \emph{Growth of the derivative of bounded analytic functions},
 Complex Variables Theory Appl. \textbf{20} (1992), no. 1--4, 221–-227.
%
\bibitem{G:BMOA} D.~Girela, \emph{Analytic functions of bounded mean
oscillation}, {Complex Function Spaces, Mekrij\"arvi 1999\/ Editor:
R.~Aulaskari. Univ. Joensuu Dept. Math. Rep. Ser. 4, Univ. Joensuu},
Joensuu (2001) pp.~61--170.
%
\bibitem{GGP} D.~Girela, C.~Gonz\'{a}lez and J.~A.~Pel\'{a}ez,
\emph{Multiplication and division by inner functions in the space of
Bloch functions}, Proc. Amer. Math. Soc. \textbf{134} (2006), no. 5,
1309–-1314.
%
\bibitem{GM1} D.~Girela and N.~Merch\'{a}n, \emph{A generalized
Hilbert operator acting on conformally invariant spaces}, Banach J.
Math. Anal. \textbf{12} (2018), no. 2, 374–-398.
%
\bibitem{GNW}
D.~Girela, M.~Nowak and P.~Waniurski, \emph{On the zeros of Bloch
functions}, Math. Proc. Cambridge Philos. Soc. \textbf{129} (2000),
no. 1, 117–-128.
%
\bibitem{GP-JFA-06} D.~Girela and J.~A.~Pel\'{a}ez, \emph{Carleson measures, multipliers and integration operators for spaces of Dirichlet
type}, J. Funct. Anal. \textbf{241} (2006), no. 1, 334– 358.
%
 \bibitem{Gu}
 V.~P.~Gurarii,
 \emph{On the factorization of absolutely convergent
 Taylor series and Fourier integrals},
 (in Russian)
 Zap. Nau\v cn. Sem. Leningrad. Otdel. Mat. Inst. Steklov. (LOMI)
 \textbf{30} (1972), 15--32.
%
\bibitem{Havin}
V.~P.~Havin, \emph{On the factorization of analytic functions smooth
up to the boundary} (Russian), Zap. Naucn. Sem. Leningrad. Otdel.
Mat. Inst. Steklov. (LOMI) \textbf{22} (1971) 202–-205.
%
\bibitem{he:kprop}
 H.~Hedenmalm,
  \emph{On the $f$- and $K$-properties of certain function spaces},
 Contemporary Math. \textbf{91} (1989), 89--91.
%
\bibitem{HKZ}
H.~Hedenmalm, B.~Korenblum \and K.~Zhu, \emph{Theory of Bergman
Spaces\/}, Graduate Texts in Mathematics, Vol. \textbf{199},
Springer, New York, Berlin, etc. 2000.
%
\bibitem{HW1}
F.~Holland and D.~Walsh, \emph{Growth estimates for functions in the
Besov spaces $A_p $}, Proc. Roy. Irish Acad. Sect. A \textbf{88}
(1988), 1--18.
%
\bibitem{kor:outer}
 B.~I.~Korenblum,
\emph{A certain extremal property of outer functions},
 (in Russian),
 Mat. Zametki \textbf{10} (1971), 53--56.
 {\sl English translation in Math. Notes \textbf{10} (1971), 456--458.}

\bibitem{kor-faivy73}
 B.~I.~Korenblum and V.~M.~Fa\u\i vy\v sevski\u\i,
 \emph{A certain class of compression operators that are
 connected with the divisibility of analytic functions},
 (in Russian),
 Ukrain. Mat.~Z. \textbf{24} (1972), 692-695.
 {\sl English translation in Ukrainian Math.~J. \textbf{24} (1973),
 559--561.}
%
\bibitem{OF} J.~M.~Ortega and J.~F\`{a}brega, \emph{Pointwise multipliers and corona type decomposition in BMOA},
 Ann. Inst. Fourier (Grenoble) \textbf{46} (1996), no. 1, 111–-137.
%
\bibitem{Pau-Pe-MZ-2009} J.~Pau and J.~A.~Pel\'{a}ez, \emph{Multipliers of M\"{o}bius invariant $Q_s$ spaces},
Math. Z. \textbf{261} (2009), no. 3, 545–-555.
%
\bibitem{Pau-Pe-2011} J.~Pau and J.~A.~Pel\'{a}ez, \emph{On the zeros of functions in Dirichlet-type spaces},
 Trans. Amer. Math. Soc. \textbf{363} (2011), no. 4, 1981–-2002.
 %
 \bibitem{Pom} Ch.~Pommerenke, \emph{Schlichte Funktionen und analytische Funktionen von beschränkter mittlerer Oszillation},
 Comment. Math. Helv. \textbf{52} (1977), no. 4, 591–-602.
%
\bibitem{RT} L.~E.~Rubel and R.~M.~Timoney,
\emph{An extremal property of the Bloch space}, Proc. Amer. Math.
Soc. \textbf{75} (1979), no.~1, 45--49.
%
\bibitem{Steg}
D.~Stegenga, \textit{Multipliers of the Dirichlet space.} Illinois
J. Math. \textbf{24} (1980), no.~1, 113--139.
%
\bibitem{Str} K.~Stroethoff,
\emph{Besov-type characterisations for the Bloch space}, Bull.
Austral. Math. Soc. \textbf{39} (1989), no. 3, 405–-420.
%
\bibitem{T} R.~M.~Timoney,
\emph{Natural function spaces}, J. London Math. Soc. (2) \textbf{41}
(1990), no. 1, 78–-88.
%
\bibitem{Vi}
S.~A.~Vinogradov, {\em{Multiplication and division in the space of
analytic functions with area integrable derivative, and in some
related spaces}} (in russian), {Zap. Nauchn. Sem. S.-Peterburg.
Otdel. Mat. Inst. Steklov. (POMI) \/} \textbf{222} (1995), {Issled.
po Linein. Oper. i Teor. Funktsii\/} \textbf{23}, 45--77, 308. {\sl
English translation in J. Math. Sci. (New York)\/ \textbf{87}, no. 5
(1997), 3806--3827.}
%
\bibitem{Wu} Z.~Wu, \textit{Carleson measures and multipliers for Dirichlet spaces.}
J. Funct. Anal. \textbf{169} (1999), no.~1, 148--163.
%
\bibitem{X}
J.~Xiao, \emph{Carleson measure, atomic decomposition and free
interpolation from Bloch space},  Ann. Acad. Sci. Fenn. Ser. A I
Math. \textbf{19} (1994), 35--46.
%
\bibitem{X-Pac}
J.~Xiao, \emph{The $Q_p$ corona theorem}, Pacific J. Math.
\textbf{194} (2000), no. 2, 491–-509.
%
\bibitem{X2}
J.~Xiao, \emph{Holomorphic $Q $ classes}, Lecture Notes in
Mathematics \textbf{1767}, Springer-Verlag, 2001.
%
\bibitem{X3}
J.~Xiao, \emph{Geometric $Q $ functions}, Frontiers in Mathematics.
Birkh\"auser Verlag, 2006.
%
\bibitem{X-2008} J.~Xiao, \emph{The $Q_p$ Carleson measure problem}, Adv. Math. \textbf{217} (2008), no. 5, 2075–-2088.
%
\bibitem{Zhao}
R.~Zhao, {\em{On logarithmic Carleson measures}}, Acta Sci. Math.
(Szeged) \textbf{69} (2003), no. 3-4, 605--618.
%
\bibitem{Z}
K.~Zhu, {\em{Analytic Besov spaces}}, J. Math. Anal. Appl.
\textbf{157} (1991), 318--336.
%
\bibitem{Zhu-MI} K.~Zhu, \emph{ A class of M\"{o}bius invariant function spaces}, Illinois J. Math.
\textbf{51} (2007), no. 3, 977-–1002.
%
\bibitem{Zhu} K.~Zhu, \emph{Operator Theory in Function Spaces\/},
Marcel Dekker, New York, 1990. Reprint: Math. Surveys and
Monographs, Vol. 138,  American Mathematical Society, Providence,
Rhode Island, 2007.
%
\bibitem{Zo} N.~Zorboska, \textit{Multiplication and Toeplitz operators on the analytic Besov
spaces.} In: \textit{More Progress in Analysis: Proceedings of the
5th International. Isaac Congress. Catania, Italy, 25 – 30 July
2005. Editors: H.~G.~W.~Begehr and F.~Nicolosi.} World Scientific
(2009), pp. 387--396.
%
\bibitem{Zy} A.~Zygmund, \emph{Trigonometric Series\/},
Vol. I and Vol. II, Second edition, Camb. Univ. Press, Cambridge,
1959.

\end{thebibliography}
\end{document}